\documentclass[a4paper,11pt]{jdg-p}
\usepackage{amssymb}
\usepackage{mathrsfs}
\usepackage{amsmath,amssymb,amsthm,latexsym,amscd,mathrsfs}
\usepackage{indentfirst}
\usepackage{stmaryrd}
 \newcommand{\ROM}[1]{\mathrm{\uppercase\expandafter{\romannumeral#1}}}
\theoremstyle{definition}

\newtheorem{thm}{Theorem}[section]
\newtheorem{lem}{Lemma}[section]

\newtheorem{rem}{Remark}[section]
\newtheorem{prop}{Proposition}[section]

\newtheorem{ack}{Acknowledgements}   
  \makeatletter
  \renewcommand{\theenumi}{\roman{enumi}}
  \renewcommand{\labelenumi}{(\theenumi)}
\makeatother
\title[Isoparametric foliation and Yau conjecture on the first eigenvalue]{\textbf{Isoparametric foliation and Yau conjecture on the first eigenvalue}}
\author[Z. Z. Tang]{Zizhou Tang}\address{School of Mathematical Sciences, Laboratory of Mathematics and Complex Systems, Beijing Normal
University, Beijing 100875, China}\email{zztang@bnu.edu.cn}
\thanks {The project is partially supported by the NSFC ( No.11071018 ) and the Program for Changjiang Scholars and Innovative
Research Team in University.}
\author[W. J. Yan]{Wenjiao Yan}
\address{School of Mathematical Sciences, Laboratory of Mathematics and Complex Systems, Beijing Normal
University, Beijing 100875, China} \email{wjyan@mail.bnu.edu.cn}

\begin{document}

\maketitle
\begin{center}
Dedicated to Professor Banghe Li on his 70th birthday.
\end{center}
\begin{abstract}
A well known conjecture of Yau states that the first eigenvalue of
every closed minimal hypersurface $M^n$ in the unit sphere
$S^{n+1}(1)$ is just its dimension $n$. The present paper shows that
Yau conjecture is true for minimal isoparametric hypersurfaces.
Moreover, the more fascinating result of this paper is that the
first eigenvalues of the focal submanifolds are equal to their
dimensions in the non-stable range.

\end{abstract}


\section{\textbf{Introduction}}
One of the most important operators acting on $C^{\infty}$ functions
on a Riemannian manifold is the Laplace-Beltrami operator. Over several decades,
research on the spectrum of the Laplace-Beltrami operator has always been a
core issue in the study of geometry.  For instance, the geometry of closed minimal submanifolds in
the unit sphere is closely related to the eigenvalue problem.

Let $(M^n,g)$ be an $n$-dimensional compact connected Riemannian manifold without boundary
and $\Delta $ be the Laplace-Beltrami operator acting on a $C^{\infty}$ function $f$ on $M$ by $\Delta f$ $=-$ div$(\nabla f)$,
the negative of divergence of the gradient $\nabla f$.
It is well known that $\Delta$ is an elliptic operator and has a discrete spectrum
$$\{0=\lambda_0(M)<\lambda_1(M)\leq \lambda_2(M)\leq \cdots \leq\lambda_k(M),\cdots, \uparrow \infty\}$$
with each eigenvalue repeated a number of times equal to its multiplicity.
As usual, we call $\lambda_1(M)$ the first eigenvalue of $M$. When $M^n$ is a minimal hypersurface
in the unit sphere $S^{n+1}(1)$, it follows from Takahashi Theorem that $\lambda_1(M)$ is
not greater than $n$.

In this connection, S.T.Yau posed in 1982 the following conjecture:
\vspace{1mm}

\noindent
\textbf{Yau conjecture (\cite{Yau}):}\,\,
{\itshape
The first eigenvalue of every closed minimal hypersurface $M^n$ in
the unit sphere $S^{n+1}(1)$ is just $n$.
}
\vspace{1mm}

The most significant breakthrough to this problem was made by Choi
and Wang (\cite{CW}). They showed that the first eigenvalue of every
(embedded) closed minimal hypersurface in $S^{n+1}(1)$ is not
smaller than $n\over 2$. As a common understanding, the calculation
of the spectrum of the Laplace-Beltrami operator, even of the first
eigenvalue, is rather complicated and difficult. Up to now, Yau
conjecture is far from being solved. In this paper, we consider a
little more restricted problem of Yau conjecture for closed minimal
isoparametric hypersurfaces $M^n$ in $S^{n+1}(1)$. As one of the
main results of this paper, we show

\begin{thm}\label{problem I}
\emph{Let $M^n$ be a closed minimal isoparametric hypersurface in $S^{n+1}(1)$. Then
$$\lambda_1(M^n)=n.$$ }
\end{thm}


Recall that a hypersurface $M^n$ in the unit sphere $S^{n+1}(1)$ is
called isoparametric if it has constant principal curvatures
(\emph{cf.} \cite{Car1}, \cite{Car2}, \cite{CR}). Let $\xi$ be a
unit normal vector field along $M^n$ in $S^{n+1}(1)$, $g$ the number
of distinct principal curvatures of $M$, $\cot \theta_{\alpha}~
(\alpha=1,...,g;~ 0<\theta_1<\cdots<\theta_{g} <\pi)$ the principal
curvatures with respect to $\xi$ and $m_{\alpha}$ the multiplicity
of $\cot \theta_{\alpha}$. Using an elegant topological method,
M\"{u}nzner proved the remarkable result that the number $g$ must be
$1, 2, 3, 4$ or $6$; $m_{\alpha}=m_{\alpha+2}$ (indices mod $g$);
$\theta_{\alpha}=\theta_1+\frac{\alpha-1}{g}\pi$ $(\alpha = 1,...,
g)$ and when $g$ is odd, $m_1=m_2$ (\emph{cf.} \cite{Mun}).

Attacking Yau conjecture, Muto-Ohnita-Urakawa (\cite{MOU}) and
Kotani (\cite{Kot}) made a breakthrough for some of the minimal
homogeneous (automatically isoparametric) hypersurfaces. More
precisely, they verified Yau conjecture for all the homogeneous
minimal hypersurfaces with $g= 1, 2, 3, 6$. However, when it came to
the case $g=4$, they were only able to deal with the cases $(m_1,
m_2)=(2,2)$ and $(1, k)$. As a matter of fact, by classification of
the homogeneous hypersurfaces with four distinct principal
curvatures, the pairs $(m_1, m_2)$ are $(1,k)$, $(2,2k-1)$,
$(4,4k-1)$, $(2,2)$, $(4,5)$ or $(6,9)$. They explained in
\cite{MOU} that ``\emph{it seems to be difficult to compute their
first eigenvalue because none of the homogeneous minimal
hypersurfaces in the unit sphere except the great sphere and the
generalized Clifford torus is symmetric or normal homogeneous}".

Furthermore, another breakthrough made by Muto (\cite{Mut}) showed that Yau conjecture is also true for some families of nonhomogeneous
minimal isoparametric hypersurfaces with four distinct principal curvatures. His remarkable result does not depend on the homogeneity
of the isoparametric hypersurfaces. However, his conclusion
covers only some isoparametric hypersurfaces with
$\min(m_1, m_2)\leq 10$. Roughly speaking, the generic families of the isoparametric hypersurfaces in the unit sphere with four distinct
principal curvatures have $\min(m_1, m_2)> 10$.

Based on all results mentioned above and the classification of isoparametric hypersurfaces in $S^{n+1}(1)$ (\emph{cf.} \cite{CCJ},
\cite{Imm}, \cite{Chi}, \cite{DN} and \cite{Miy}),
we show our Theorem \ref{problem I} by establishing the following

\begin{thm}\label{thm1 hypersurface}
\emph{Let $M^n$ be a closed minimal isoparametric hypersurface in the unit sphere $S^{n+1}(1)$ with four
distinct principal curvatures and $m_1, m_2\geq 2$. Then
$$\lambda_1(M^n)=n.$$ }
\end{thm}

\begin{rem}\label{rem 1}
For isoparametric hypersurfaces in the unit spheres with $g = 1, 2, 3$,
Cartan classified them to be homogeneous (\emph{cf. }\cite{Car1}, \cite{Car2}); for $g = 6$, Dorfmeister- Neher (\cite{DN})
and Miyaoka (\cite{Miy}) showed that they are homogeneous. Thus the results of \cite{MOU} and \cite{Kot} complete the proof of
Theorem \ref{problem I} in cases $g=1,2,3,6$. Moreover, Takagi (\cite{Tak1}) asserted that the isoparametric hypersurface
with $g=4$ and multiplicities $(1,k)$ must be homogeneous. By virtue of \cite{MOU}, Theorem \ref{problem I} is true for the
case $(1,k)$. Therefore, Theorem \ref{thm1 hypersurface} completes in a direct way the proof of Theorem \ref{problem I}.
\end{rem}

\begin{rem}
For isoparametric hypersurfaces with $g=4$, Cecil-Chi-Jensen
(\cite{CCJ}), Immervoll (\cite{Imm}) and Chi (\cite{Chi}) proved a
far reaching result that they are either homogeneous or of
OT-FKM-type except possibly for the case $(m_1, m_2)=(7,8)$.
Actually, Theorem \ref{thm1 hypersurface} depends only on the values
of $(m_1, m_2)$, but not on the homogeneity. Besides, our method is
also applicable to the case $g=6$.
\end{rem}

\begin{rem}\label{Chern conjecture}
Chern conjectured that a closed, minimally immersed hypersurface in
$S^{n+1}(1)$, whose second fundamental form has constant length, is
isoparametric (\emph{cf.} \cite{GT}). If this conjecture is proven,
we would have settled Yau conjecture for the minimal hypersurface
whose second fundamental form has constant length, which gives us
more confidence in Yau conjecture.
\end{rem}

The more fascinating part of this paper is the determination of the
first eigenvalues of the focal submanifolds in $S^{n+1}(1)$, which
relays on the deeper geometric properties of the isoparametric
foliation.

To state our Theorem \ref{thm2 focal submanifold} clearly, let us make some preliminaries.
A well known result of Cartan states that isoparametric hypersurfaces come as a family of
parallel hypersurfaces. To be more specific, given an isoparametric hypersurface $M^n$ in $S^{n+1}(1)$ and a smooth field $\xi$ of unit
normals to $M$, for each $x\in M$ and $\theta\in \mathbb{R}$, we can define
$\phi_{\theta}: M^n\rightarrow S^{n+1}(1)$ by
$$\phi_{\theta}(x)=\cos \theta~ x +\sin \theta~ \xi(x).$$
Clearly, $\phi_{\theta}(x)$ is the point at an oriented distance
$\theta$ to $M$ along the normal geodesic through $x$. If
$\theta\neq \theta_{\alpha}$ for any $\alpha=1,...,g$,
$\phi_{\theta}$ is a parallel hypersurface to $M$ at an oriented
distance $\theta$, which we will denote by $M_{\theta}$
henceforward. If $\theta= \theta_{\alpha}$ for some
$\alpha=1,...,g$, it is easy to find that for any vector $X$ in the
principal distributions $E_{\alpha}(x)=\{X\in T_xM ~|~A_{\xi}X=\cot
\theta_{\alpha}X\}$, where $A_{\xi}$ is the shape operator with
respect to $\xi$, $(\phi_{\theta})_{\ast}X=(\cos \theta-\sin \theta
\cot
\theta_{\alpha})X=\frac{sin(\theta_{\alpha}-\theta)}{sin\theta_{\alpha}}
X=0$. In other words, in case $\cot \theta=\cot \theta_{\alpha}$ is
a principal curvature of $M$, $\phi_{\theta}$ is not an immersion,
which is actually a \emph{focal submanifold} of codimension
$m_{\alpha}+1$ in $S^{n+1}(1)$.

M\"{u}nzner asserted that regardless of the number
of distinct principal curvatures of $M$, there are only two distinct
focal submanifolds in a parallel family of isoparametric hypersurfaces,
and every isoparametric hypersurface is a tube of constant radius over each focal submanifold.
Denote by $M_1$ the focal submanifold in $S^{n+1}(1)$ at an oriented distance $\theta_1$
along $\xi$ from $M$ with codimension $m_1+1$, $M_2$ the focal submanifold in $S^{n+1}(1)$ at
an oriented distance $\frac{\pi}{g}-\theta_1$ along $-\xi$ from $M$ with codimension $m_2+1$.
In virtue of Cartan's identity, one sees that
the focal submanifolds $M_1$ and $M_2$ are minimal in $S^{n+1}(1)$ (\emph{cf.} \cite{CR}).

Another main result of the present paper concerning the first
eigenvalues of focal submanifolds in the non-stable range
(\emph{cf}. \cite{HH}), is stated as follows.

\begin{thm}\label{thm2 focal submanifold}
\emph{Let $M_1$ be the focal submanifold of an isoparametric
hypersurface with four distinct principal curvatures in the unit
sphere $S^{n+1}(1)$ with codimension $m_1+1$. If $\dim M_1\geq
\frac{2}{3}n+1$, then
$$\lambda_1(M_1)=\dim M_1$$ with multiplicity $n+2$.
A similar conclusion holds for $M_2$ under an analogous condition.}
\end{thm}

As a simple application of Theorem \ref{thm2 focal submanifold}, we
obtain that each focal submanifold of isoparametric hypersurfaces
with $g=4$, $(m_1, m_2)=(4, 5)$ or $(7,8)$ has its dimension as the
first eigenvalue.


Bearing the results above in mind, in conjunction with the
classification results of \cite{CCJ} \cite{Chi} which stated that
except for the case $(m_1, m_2)=(7, 8)$, the isoparametric
hypersurfaces in $S^{n+1}(1)$ with four distinct principal
curvatures are either homogeneous with $(m_1, m_2)=(2,2), (4,5)$ or
of OT-FKM-type, we will look into the focal submanifolds of
OT-FKM-type and give their first eigenvalues.

We now recall the construction of the isoparametric hypersurfaces of
OT-FKM-type. For a symmetric Clifford system $\{P_0,\cdots,P_m\}$ on
$\mathbb{R}^{2l}$, \emph{i.e.}, $P_i$'s are symmetric matrices
satisfying $P_iP_j+P_jP_i=2\delta_{ij}I_{2l}$, Ferus, Karcher and
M\"{u}nzner (\cite{FKM}) constructed a polynomial $F$ on
$\mathbb{R}^{2l}$:
\begin{eqnarray}\label{FKM isop. poly.}
&&\qquad F:\quad \mathbb{R}^{2l}\rightarrow \mathbb{R}\nonumber\\
&&F(x) = |x|^4 - 2\displaystyle\sum_{i = 0}^{m}{\langle
P_ix,x\rangle^2}.
\end{eqnarray}

It turns out that each level hypersurface of $f=F|_{S^{2l-1}}$, \emph{i.e.},
the preimage of some regular value of $f$, has four distinct
constant principal curvatures. Choosing $\xi=\frac{\nabla f}{|\nabla f|}$, we find $M_1=f^{-1}(1)$, $M_2=f^{-1}(-1)$,
which have codimensions $m_1+1$ and $m_2+1$ in $S^{n+1}(1)$, respectively. The multiplicity pairs $(m_1, m_2)$
of the OT-FKM-type are $(m, l-m-1)$, provided
$m>0$ and $l-m-1> 0$, where $l = k\delta(m)$ $(k=1,2,3,...)$ and $\delta(m)$ is the dimension
of an irreducible module of the Clifford algebra $C_{m-1}$ .
In the following, we list the values of $\delta(m)$ corresponding to $m$:

\begin{center}
\begin{tabular}{|c|c|c|c|c|c|c|c|c|c|}
\hline
$m$ & 1 & 2 & 3 & 4 & 5 & 6 & 7 & 8 & $\cdots$ $m$+8 \\
\hline
$\delta(m)$ & 1 & 2 & 4 & 4 & 8 & 8 & 8 & 8 & ~16$\delta(m)$\\
\hline
\end{tabular}
\end{center}
\vspace{2mm}

Firstly, we focus on the focal submanifold $M_2$. If $3\dim M_2\geq
2n+3$, or equivalently, $m_1\geq\frac{1}{2}(m_2+3)$, Theorem
\ref{thm2 focal submanifold} gives $\lambda_1(M_2)=\dim
M_2=2m_1+m_2$. The assumption $3\dim M_2\geq 2n+3$ is essential. For
instance, Solomon (\cite{Sol}) constructed an eigenfunction on the
focal submanifold $M_2$ of OT-FKM-type, which has $4m$ as an
eigenvalue. It follows that $\lambda_1(M_2)\leq 4m$. Therefore, in
the stable range $3\dim M_2< 2(n+1)-2$, \emph{i.e.} $m_1<
\frac{1}{2}m_2$, $\lambda_1(M_2)< 2m_1+m_2=\dim M_2.$ Only three
cases are left to estimate: $m_1=\frac{1}{2}m_2$,
$m_1=\frac{1}{2}(m_2+1)$ and $m_1=\frac{1}{2}(m_2+2)$, which are
actually $(m_1, m_2)=(1,1)$, $(1, 2)$, $(2, 3)$, $(3, 4)$, $(4, 7)$,
$(5,10)$ and $(8, 15)$.

Next, we will be concerned with the focal submanifold $M_1$. Fortunately, the condition in
Theorem \ref{thm2 focal submanifold} is almost satisfied. Equivalently, the first eigenvalue of the focal submanifold $M_1$ of
those OT-FKM-type can be determined completely.
By analyzing the conditions $m_1\geq 1$, $m_2\geq 1$ and $m_2<\frac{1}{2}(m_1+3)$,
we find that there are only five cases left, that is, $(m_1, m_2)=(1, 1)$, $(2, 1)$, $(4, 3)$, $(5, 2)$ and $(6,1)$.
In view of \cite{FKM}, the families for multiplicities $(2, 1)$, $(5, 2)$, $(6,1)$ and one of the $(4, 3)$-families are congruent
to those with multiplicities $(1, 2)$, $(2, 5)$, $(1, 6)$ and $(3, 4)$, respectively, and the focal submanifolds interchange.
For the case $(2, 5)$, an effective estimate can be given by \cite{Sol}, while for the cases $(1, 2)$ and $(1, 6)$, the following
proposition determines the first eigenvalues.

\begin{prop}\label{M1 for other}
\emph{Let $M_2$ be the focal submanifold of OT-FKM-type defined before with $(m_1, m_2)=(1, k)$. The following equality is valid
\begin{equation*}
\lambda_1(M_2)=\min\{4, 2+k\}.
\end{equation*}}
\end{prop}

As mentioned before, Takagi (\cite{Tak1}) asserted that the
isoparametric hypersurface with $g=4$ and multiplicity $(1,k)$ must
be homogeneous. Thus the corresponding focal submanifold of
isoparametric hypersurface with four distinct principal curvatures
and $\min\{m_1, m_2\}=1$ has $\min\{4, 2+k\}$ as its first
eigenvalue.


\vspace{2mm}

At last, we would like to propose a problem on the first eigenvalue
of the minimal submanifolds with dimensions in the non-stable range
in $S^{n+1}(1)$, which could be regarded as an extension of Yau
conjecture. \vspace{2mm}

\noindent \textbf{Problem:}\,\, {\itshape Let $M^d$ be a closed
minimal submanifold in the unit sphere $S^{n+1}(1)$ with
$d\geq \frac{2}{3}n+1$.
Is it true that
$$\lambda_1(M^d)=d~?$$
}


\section{\textbf{The first eigenvalue of the minimal isoparametric hypersurface}}

Let $\phi: M^n \rightarrow S^{n+1}(1)(\subset \mathbb{R}^{n+2})$ be a closed isoparametric hypersurface with $g$ distinct
principal curvatures in $S^{n+1}(1)$ and
$\xi$ be a smooth field of unit normals to $M$. Again, denote by $E_{\alpha}~ (\alpha~=1,..., g)$ the principal
distribution on $M$, \emph{i.e.}, the eigenspace of the shape operator $A_{\xi}$ corresponding to the eigenvalue
$\cot \theta_{\alpha}~ (0<\theta_1<...<\theta_{g} <\pi)$. The parallel hypersurface $M_{\theta}$ at an oriented distance
$\theta$ from $\phi$ is defined by $\phi_{\theta}: M^n \rightarrow S^{n+1}(1)~(-\pi<\theta<\pi, \cot\theta\neq \cot \theta_{\alpha})$,
$$\phi_{\theta}(x) = \cos \theta~ x + \sin \theta~ \xi(x).$$
At first, let us prepare some formulae:

For $X\in E_{\alpha}$, it is easy to see
\begin{equation}\label{phi theta star}
(\phi_{\theta})_{\ast} X=\frac{sin(\theta_{\alpha}-\theta)}{sin
\theta_{\alpha}} \widetilde{X},
\end{equation}
where $\widetilde{X} \sslash X$ as vectors in $\mathbb{R}^{n+2}$.


Let $H$ be the mean curvature of $M^n$ in $S^{n+1}(1)$ with respect to $\xi$. Clearly,
\begin{eqnarray}\label{mean curvature}
nH&=&\sum_{\alpha=1}^{g}m_{\alpha}\cot\theta_{\alpha}\\
&=&\left\{ \begin{aligned}
m_1g\cot(g\theta_1)\qquad\qquad\qquad\qquad for~ g ~odd\\
\frac{m_1g}{2}\cot\frac{g\theta_1}{2}-\frac{m_2g}{2}\tan\frac{g\theta_1}{2}\qquad\qquad for~ g~even
\end{aligned}\right.
\nonumber
\end{eqnarray}

In order to estimate the eigenvalues of $M$, we would recall a theorem that will play a crucial role in our work as Muto did in \cite{Mut}.
\vspace{1mm}

\noindent \textbf{Theorem (Chavel and Feldman \cite{CF}, Ozawa
\cite{Oza})}\,\, {\itshape Let $V$ be a closed, connected smooth
Riemannian manifold and $W$ a closed submanifold of $V$. For any
sufficiently small $\varepsilon>0$, set $W(\varepsilon)=\{x\in V:~
dist(x, W)<\varepsilon\}$. Let $\lambda^D_k(\varepsilon)$ $(k=1, 2,
...)$ be the $k$-th eigenvalue of the Laplace-Beltrami operator on
$V-W(\varepsilon)$ under the Dirichlet boundary condition. If $\dim
V\geq \dim W+ 2$, then for any $k = 1, 2$...
\begin{equation}\label{lambda k limit}
\lim_{\varepsilon\to 0}\lambda^D_k(\varepsilon)=\lambda_{k-1}(V).
\end{equation}}

We will apply this theorem to the case $V=S^{n+1}(1)$ and $W=M_1\cup M_2$, the union of the focal submanifolds.
By estimating the eigenvalue
$\lambda_{k}(M^n)$ from below, we can prove Theorem \ref{thm1 hypersurface}.
\vspace{3mm}

\noindent \textbf{Theorem 1.2.}\, \emph{Let $M^n$ be a closed
minimal isoparametric hypersurface in the unit sphere $S^{n+1}(1)$
with four distinct principal curvatures and $m_1, m_2\geq 2$. Then
$$\lambda_1(M^n)=n.$$ }
\vspace{3mm}

\noindent
\textbf{\emph{Proof}}. For sufficiently small $\varepsilon> 0$, set
$$M(\varepsilon)=\bigcup_{\theta\in[-\frac{\pi}{4}+\theta_1+\varepsilon, ~\theta_1-\varepsilon]}M_{\theta}.$$
Clearly, $M(\varepsilon)$ is a domain of $S^{n+1}(1)$ obtained by excluding $\varepsilon$-neighborhoods of
$M_1$ and $M_2$ from $S^{n+1}(1)$. Alternatively, it can also be regarded as a tube around the minimal isoparametric hypersurface $M$.
According to the theorem of Chavel, Feldman and Ozawa,
\begin{equation}\label{Chavel}
\lim_{\varepsilon\to 0}\lambda^D_{k+1}(M(\varepsilon))=\lambda_{k}(S^{n+1}(1)),
\end{equation} we
need to estimate $\lambda^D_{k+1}(M(\varepsilon))$ from above in terms of $\lambda_{k}(M^n)$.

Let $\Big\{\widetilde{e}_{\alpha,i}~ \boldsymbol{|} ~ i=1,...,m_{\alpha},~ \alpha=1,..,4,~ \widetilde{e}_{\alpha, i}\in E_{\alpha}\Big\}$
be a local orthonormal frame field on $M$.
Then $$\Big\{\frac{\partial}{\partial\theta},~ e_{\alpha,i}~\boldsymbol{|}~ e_{\alpha,i}=\frac{\sin\theta_{\alpha}}{\sin(\theta_{\alpha}-\theta)}\widetilde{e}_{\alpha,i},~
i=1,...,m_{\alpha},~ \alpha=1,..,4,~ \theta\in[-\frac{\pi}{4}+\theta_1+\varepsilon, ~\theta_1-\varepsilon]\Big\}$$ is a local orthonormal frame field on $M(\varepsilon)$.
From the formula (\ref{phi theta star}), we derive immediately that the volume element of $M(\varepsilon)$
can be expressed in terms of the volume element of $M$:
\begin{equation}\label{volume}
dM(\varepsilon)=\frac{\sin^{m_1}2(\theta_1-\theta)\cos^{m_2}2(\theta_1-\theta)}{\sin^{m_1}2\theta_1\cos^{m_2}2\theta_1}d\theta dM
\end{equation}

Following \cite{Mut}, let $h$ be a nonnegative, increasing smooth function on $[0, \infty)$ satisfying $h=1$ on
$[2, \infty)$ and $h=0$ on $[0, 1]$. For sufficiently small $\eta>0$, let $\psi_{\eta}$ be a nonnegative
smooth function on $[\eta, \frac{\pi}{2}-\eta]$ such that

$(i)$ $\psi_{\eta}(\eta)=\psi_{\eta}(\frac{\pi}{2}-\eta)=0$,

$(ii)$ $\psi_{\eta}$ is symmetric with respect to $x=\frac{\pi}{4}$

$(iii)$ $\psi_{\eta}(x)=h(\frac{x}{\eta})$ on $[\eta, \frac{\pi}{4}]$.

\noindent
Let $f_k$ $(k=0, 1, ...)$ be the $k$-th eigenfunctions on $M$ which are orthogonal to
each other with respect to the square integral inner product on $M$ and $L_{k+1} =Span\{f_0, f_1,..., f_k\}$.

For each fixed $\theta\in[-\frac{\pi}{4}+\theta_1+\varepsilon, ~\theta_1-\varepsilon]$,
denote $\pi=\pi_{\theta}=\phi^{-1}_{\theta}: M_{\theta}\rightarrow M$. Then any
$\varphi\in L_{k+1}$ on $M$ can give rise to
a function $\Phi_{\varepsilon} : M(\varepsilon)\rightarrow \mathbb{R}$ by
$$\Phi_{\varepsilon}(x) = \psi_{2\varepsilon}(2(\theta_1-\theta))(\varphi\circ\pi)(x),$$
where $\theta$ is characterized by $x\in M_{\theta}$, $\theta\in[-\frac{\pi}{4}+\theta_1+\varepsilon, ~\theta_1-\varepsilon]$.
It is evident to see that $\Phi_{\varepsilon}$ is a smooth function on $M(\varepsilon)$ satisfying the Dirichlet boundary condition
and square integrable.

By the mini-max principle, we have:
\begin{equation}\label{min max}
\lambda^D_{k+1}(M(\varepsilon))\leq \sup_{\varphi\in L_{k+1}}\frac{\|\nabla \Phi_{\varepsilon}\|_2^2}{\|\Phi_{\varepsilon}\|_2^2}.
\end{equation}
In the following, we will concentrate on the calculation of $\frac{\|\nabla \Phi_{\varepsilon}\|_2^2}{\|\Phi_{\varepsilon}\|_2^2}$.
Observing that the normal geodesic starting from $M$ is perpendicular to each parallel hypersurface $M_{\theta}$, we obtain
$$\|\nabla \Phi_{\varepsilon}\|_2^2=\int_{M(\varepsilon)}4(\psi^{\prime}_{2\varepsilon})^2\varphi(\pi)^2dM(\varepsilon)
+\int_{M(\varepsilon)}\psi^{2}_{2\varepsilon}|\nabla \varphi(\pi)|^2dM(\varepsilon).$$
On the other hand, a simple calculation leads to
\begin{eqnarray}\label{Phi length}
\|\Phi_{\varepsilon}\|_2^2&=&\int_{M(\varepsilon)}\psi^2_{2\varepsilon}(2(\theta_1-\theta))\varphi(\pi(x))^2 dM(\varepsilon)\nonumber\\
&=&\int_M\int_{-\frac{\pi}{4}+\theta_1+\varepsilon}^{\theta_1-\varepsilon}\psi^2_{2\varepsilon}(2(\theta_1-\theta))
\frac{\sin^{m_1}2(\theta_1-\theta)\cos^{m_2}2(\theta_1-\theta)}{\sin^{m_1}2\theta_1\cos^{m_2}2\theta_1}\varphi(\pi(x))^2d\theta dM\nonumber\\
&=&\frac{\|\varphi\|^2_2}{2\sin^{m_1}2\theta_1\cos^{m_2}2\theta_1}\Big(\int_{2\varepsilon}^{\frac{\pi}{2}-2\varepsilon}\psi^2_{2\varepsilon}(x)\sin^{m_1}x\cos^{m_2}x~dx\Big).\nonumber
\end{eqnarray}
\noindent
For the sake of convenience, let us decompose
\begin{equation}\label{Phi}
\frac{\|\nabla \Phi_{\varepsilon}\|_2^2}{\|\Phi_{\varepsilon}\|_2^2}=I(\varepsilon)+II(\varepsilon),
\end{equation}
with
\begin{eqnarray}\label{I}
I(\varepsilon)
&=& \frac{\int_{M(\varepsilon)}4(\psi^{\prime}_{2\varepsilon})^2\varphi(\pi)^2dM(\varepsilon)}{\int_{M(\varepsilon)}(\psi_{2\varepsilon})^2\varphi(\pi)^2~dM(\varepsilon)}\\
&=&\frac{4\int_{2\varepsilon}^{\frac{\pi}{2}-2\varepsilon}(\psi^{\prime}_{2\varepsilon}(x))^2\sin^{m_1}x\cos^{m_2}x~dx}
{\int_{2\varepsilon}^{\frac{\pi}{2}-2\varepsilon}\psi^2_{2\varepsilon}(x)\sin^{m_1}x\cos^{m_2}x~ dx}\nonumber
\end{eqnarray}
and
\begin{equation}\label{II}
II(\varepsilon)
 =\frac{\int_{M(\varepsilon)}\psi^{2}_{2\varepsilon}|\nabla \varphi(\pi)|^2dM(\varepsilon)}
 {\int_{M(\varepsilon)}\psi^2_{2\varepsilon}\varphi(\pi)^2 dM(\varepsilon)}.
\end{equation}

We shall take the first step by claiming that
\begin{equation}\label{limit I}
\lim_{\varepsilon\rightarrow 0}I(\varepsilon)=0.
\end{equation}
In fact, for the smooth function $h$, we have a positive number $C$
such that $|h^{\prime}|\leq C$. It follows immediately that
$|\psi_{\eta}^{\prime}(x)|=|\frac{1}{\eta}h^{\prime}(\frac{x}{\eta})|\leq
\frac{1}{\eta}C$ for $x\in[\eta, \frac{\pi}{4}]$. Under the
assumption $\min\{m_1, m_2\}\geq 2$, we deduce
\begin{eqnarray*}\label{limit I proof}
&&\int_{2\varepsilon}^{\frac{\pi}{2}-2\varepsilon}(\psi^{\prime}_{2\varepsilon}(x))^2\sin^{m_1}x\cos^{m_2}x~dx\\
&\leq& \int_{2\varepsilon}^{4\varepsilon}(\psi^{\prime}_{2\varepsilon}(x))^2\sin^{2}x~dx + \int_{\frac{\pi}{2}-4\varepsilon}^{\frac{\pi}{2}-2\varepsilon}(\psi^{\prime}_{2\varepsilon}(x))^2\cos^{2}x~dx\nonumber\\
&\leq&
\frac{C^2}{4}\int_{2\varepsilon}^{4\varepsilon}\frac{\sin^{2}x}{{\varepsilon}^2}~dx
+
\frac{C^2}{4}\int_{\frac{\pi}{2}-4\varepsilon}^{\frac{\pi}{2}-2\varepsilon}\frac{\cos^{2}x}{{\varepsilon}^2}~dx,\nonumber
\end{eqnarray*}
from which it follows that the numerator of $I(\varepsilon)$ in
(\ref{I}) approaches to $0$ as $\varepsilon$ goes to $0$. On the
other hand, the denominator of $I(\varepsilon)$ approaches to a
non-zero number as $\varepsilon$ goes to $0$. Thus the claim
(\ref{limit I}) is established. \vspace{3mm}

Next, we turn to the estimation of $II(\varepsilon)$.

Decompose $\nabla \varphi = Z_1 + Z_2 + Z_3 + Z_4 \in E_1\oplus E_2\oplus E_3\oplus E_4$, and
set $k_{\alpha}=\frac{sin(\theta_{\alpha}-\theta)}{sin \theta_{\alpha}}$ for $\alpha=1,...,4$.
Using the following identity
\begin{equation}\label{inner product}
\langle \nabla \varphi(\pi), X \rangle = \langle \nabla \varphi, \pi_{\ast}X\rangle,\quad for~ any ~ X\in T_xM_{\theta},
\end{equation}
we have
\begin{equation}\label{nabla varphi(pi)}
\left\{ \begin{aligned}
\quad|\nabla \varphi|^2~~~&=|Z_1|^2+|Z_2|^2+|Z_3|^2+|Z_4|^2\\
|\nabla \varphi(\pi)|^2&=\frac{1}{k_1^2}|Z_1|^2+\frac{1}{k_2^2}|Z_2|^2+\frac{1}{k_3^2}|Z_3|^2+\frac{1}{k_4^2}|Z_4|^2.
\end{aligned}\right.
\end{equation}
Moreover, for simplicity, for $\alpha=1,...,4$, define
\begin{eqnarray}\label{K alpha}
K_{\alpha}&:=&\int_{-\frac{\pi}{4}+\theta_1}^{\theta_1}\frac{\sin^{m_1}2(\theta_1-\theta)\cos^{m_2}2(\theta_1-\theta)}{k^2_{\alpha}}d\theta\\
&=&\sin^2\theta_{\alpha}\int_0^{\frac{\pi}{4}}\frac{\sin^{m_1}2x\cos^{m_2}2x}{\sin^2(\frac{\alpha-1}{4}\pi+x)} dx,\nonumber\\
G&:=&\int_0^{\frac{\pi}{2}}\sin^{m_1}x\cos^{m_2}x~dx\\
&=&2\int_0^{\frac{\pi}{4}}\sin^{m_1}2x\cos^{m_2}2x~dx.\nonumber
\end{eqnarray}
Let $\displaystyle K=\max_{\alpha}\{K_{\alpha}\}$. Then combining with (\ref{Phi}), (\ref{I}), (\ref{II}), (\ref{limit I}), (\ref{nabla varphi(pi)}), (\ref{K alpha}) and (15), we arrive at
\begin{equation}\label{limit II}
\lim_{\varepsilon\rightarrow 0}\frac{\|\nabla \Phi_{\varepsilon}\|_2^2}{\|\Phi_{\varepsilon}\|_2^2}=
\frac{ \sum_{\alpha}K_{\alpha}\|Z_{\alpha}\|_2^2}{\|\varphi\|_2^2\cdot \frac{1}{2}G}
\leq \frac{2K}{G}\cdot \frac{\|\nabla \varphi\|_2^2}{\|\varphi\|_2^2}
\end{equation}

\noindent
Therefore, putting (\ref{Chavel}), (\ref{min max}) and (\ref{limit II})together, we see that
\begin{equation}\label{sphere M}
\lambda_k(S^{n+1}(1))
=\lim_{\varepsilon\rightarrow 0}\lambda_{k+1}^D(M(\varepsilon))
\leq\lim_{\varepsilon\rightarrow 0}\sup_{\varphi\in L_{k+1}}\frac{\|\nabla \Phi_{\varepsilon}\|_2^2}{\|\Phi_{\varepsilon}\|_2^2}
\leq\lambda_k(M^n) \frac{2K}{G}.
\end{equation}

\noindent Comparing the leftmost side with the rightmost side of
(\ref{sphere M}), it is sufficient to complete the proof of Theorem
\ref{thm1 hypersurface}, if we can verify the inequality
\begin{equation}\label{KG}
K<\frac{n+2}{n}G.
\end{equation}
Since then, $\lambda_{n+3}(S^{n+1}(1))=2(n+2)<
\lambda_{n+3}(M^n)\cdot\frac{2(n+2)}{n}$, which implies immediately
that $\lambda_{n+3}(M^n)> n$. Recall that $n$ is an eigenvalue of
$M^n$ with multiplicity at least $n+2$. Therefore, the first
eigenvalue of $M^n$ must be $n$ with multiplicity $n+2$.

We are now in a position to verify the inequality
(\ref{KG}), which is equivalent to
\begin{equation}\label{K alpha <G}
K_{\alpha}<\frac{n+2}{n}G, \quad for~ each~ \alpha=1,2,3,4.
\end{equation}
\noindent First, we observe that the certifications for $K_2$ and
$K_3$ are similar, so are for $K_1$ and $K_4$. Thus we just need to
give two verifications.

$(i)$ ~Given $0<x<\frac{\pi}{4}$, since $0<\theta_1<\frac{\pi}{4}$,
it follows straightforwardly that
$$K_2<2\sin^2\theta_2\int_0^{\frac{\pi}{4}}\sin^{m_1}2x\cos^{m_2}2x~dx<2\int_0^{\frac{\pi}{4}}\sin^{m_1}2x\cos^{m_2}2x~dx=G.$$
Similarly, we have
$$K_3<2\sin^2\theta_3\int_0^{\frac{\pi}{4}}\sin^{m_1}2x\cos^{m_2}2x~dx<2\int_0^{\frac{\pi}{4}}\sin^{m_1}2x\cos^{m_2}2x~dx=G.$$

$(ii)$~Express $K_1$ and $G$ in terms of the Beta function $B(x, y)= 2\int_{0}^{\frac{\pi}{2}}\sin^x\theta\cos^y\theta~d\theta$:
\begin{eqnarray}\label{K1}
K_1&=&\sin^2\theta_1\int_0^{\frac{\pi}{4}}\frac{\sin^{m_1}2x\cos^{m_2}2x}{\sin^2x} dx\\
&=& \frac{1}{2}\sin^2\theta_1\Big[B(\frac{m_1-1}{2}, \frac{m_2+1}{2})+ B(\frac{m_1-1}{2}, \frac{m_2+2}{2})\Big],\nonumber
\end{eqnarray}
\begin{equation}\label{G}
G~~=\int_0^{\frac{\pi}{2}}\sin^{m_1}x\cos^{m_2}x~dx=\frac{1}{2}B(\frac{m_1+1}{2}, \frac{m_2+1}{2}).\qquad\quad
\end{equation}
Using the properties of Beta function and Gamma function:
$$B(x,y)=\frac{\Gamma(x)\Gamma(y)}{\Gamma(x+y)}\quad and \quad\Gamma(x+1)=x\Gamma(x)\quad for ~any ~ x>0,$$
it follows from (\ref{K1}) and (\ref{G}) that
$$\frac{K_1}{G}=\sin^2\theta_1\cdot\frac{m_1+m_2}{m_1-1}\cdot
\Big(1+\frac{\Gamma(\frac{m_2+2}{2})\Gamma(\frac{m_1+m_2}{2})}{\Gamma(\frac{m_2+1}{2})\Gamma(\frac{m_1+m_2+1}{2})}\Big)$$
Define $$S(m_1,m_2):=\frac{\Gamma(\frac{m_2+2}{2})\Gamma(\frac{m_1+m_2}{2})}{\Gamma(\frac{m_2+1}{2})\Gamma(\frac{m_1+m_2+1}{2})}$$
and $$A(m_1,m_2):=\frac{n+2}{n}\frac{1}{\sin^2\theta_1}\frac{m_1-1}{m_1+m_2}.$$ Then it is clear that
\begin{equation}\label{SA}
K_1<\frac{n+2}{n}G~\Longleftrightarrow~1+S(m_1,m_2)<A(m_1,m_2).
\end{equation}

We conclude this section with establishing two inequalities $S(m_1,m_2)<1$ and $A(m_1,m_2)\geq 2$.

\begin{lem}\label{S}

\emph{The multiplicities $m_1, m_2$ of the principal curvatures of
isoparametric hypersurfaces with four distinct principal curvatures
with $m_1, m_2 \geq 2$ satisfy
$$S(m_1,m_2)<1.$$}
\end{lem}

\noindent \emph{\textbf{Proof}}. Recall a well known result that
when $g=4$, $m_1$ and $m_2$ can not be both even except for $(2,2)$
(\emph{cf.} \cite{Mun}, \cite{Abr}, \cite{Tan}). It suffices to
estimate $S(m_1,m_2)$ in the following three cases.

\noindent
Case 1: When $(m_1, m_2)=(2,2),$
$$S(2, 2):=\frac{\Gamma(2)\Gamma(2)}{\Gamma(\frac{3}{2})\Gamma(\frac{5}{2})}=\frac{8}{3\pi}<1.$$

\noindent
Case 2: When $m_1=2p+1$, it is obvious that
$$S(m_1,m_2)=\frac{\frac{m_2+1}{2}\cdot(\frac{m_2+1}{2}+1)\cdots (\frac{m_2+1}{2}+p-1)}{\frac{m_2+2}{2}\cdot(\frac{m_2+2}{2}+1)\cdots (\frac{m_2+2}{2}+p-1)}<1;$$

\noindent
Case 3: When $m_1=2p$, $m_2=2q+1$, for simplicity, we define
$$T(p,q):=S(m_1,m_2)=\frac{(2q+1)!!(2p+2q-1)!!\cdot \pi}{q!(p+q)!\cdot 2^{p+2q+1}}.$$
It is straightforward to see that $T(p,q)$ is strictly decreasing with $p$ for a fixed $q$, and strictly increasing with $q$ for a fixed $p$.
It follows that
$$T(p,q)<T(p-1,q)<\cdots <T(1,q)<T(1,q+1)<\cdots<T(1,\infty).$$
Using the Stirling Formula:
$$\lim_{n\rightarrow \infty}\frac{n!}{\sqrt{2\pi n}(\frac{n}{e})^n}=1,$$
we obtain that
\begin{eqnarray}\label{T(1,q)}
T(1,\infty)
&=&\lim_{q\rightarrow\infty}\frac{[(2q+1)!]^2\pi}{(q!)^3(q+1)!2^{4q+2}}\nonumber\\
&=&\lim_{q\rightarrow\infty}\frac{(2q+1)^{3q+\frac{3}{2}}}{(2q)^{3q+\frac{3}{2}}}\cdot \frac{(2q+1)^{q+\frac{3}{2}}}{(2q+2)^{q+\frac{3}{2}}}\cdot \frac{1}{e}\nonumber\\
&=& e^{\frac{3}{2}}\cdot \frac{1}{e^{\frac{1}{2}}}\cdot\frac{1}{e}\nonumber\\
&=&1.\nonumber
\end{eqnarray}
This completes the proof of Lemma \ref{S}.\hfill $\Box$
\vspace{3mm}

Lemma \ref{S} reduces the proof of (\ref{K alpha <G}) for $K_1$ to proving that $A(m_1,m_2)\geq 2$.

Since $M^n$ is the minimal isoparametric hypersurface in $S^{n+1}(1)$, from Formula (\ref{mean curvature}),
we derive that $\sin^2\theta_1=\frac{1}{2}(1-\frac{\sqrt{m_2}}{\sqrt{m_1+m_2}})$.
On the other hand, in our case $g=4$, we have
$n=\frac{g}{2}(m_1+m_2)=2(m_1+m_2)$, thus
$$A(m_1,m_2)=\frac{m_1-1}{m_1+m_2}\cdot\frac{m_1+m_2+1}{m_1+m_2}\cdot\frac{2}{1-\frac{\sqrt{m_2}}{\sqrt{m_1+m_2}}}.$$
A simple calculation shows
\begin{equation*}\label{A geq 2}
A(m_1,m_2)\geq 2\Longleftrightarrow m_2(m_1+m_2)^3\geq
(m_2^2+m_1m_2+m_2+1)^2.
\end{equation*}
It is not difficult to see that the following three inequalities
guarantee the right hand of the equivalence above.
\begin{equation*}
\left\{ \begin{array}{lll}
~3m_1\geq 2m_1+2\\
3{m_1}^2\geq m_1^2+2m_1+3\\
\quad m_1^3\geq 2m_1+3
\end{array}\right.
\end{equation*}

Fortunately, the last three inequalities are satisfied
simultaneously if $m_1\geq 2$. Thus $1+S(m_1,m_2)<2\leq A(m_1,m_2)$
under the assumption $\min\{m_1, m_2\}\geq 2$, equivalently, the
inequality $K_1<\frac{n+2}{n}G$ we required holds true.

Similarly, $K_4<\frac{n+2}{n}G$.

The proof of Theorem \ref{thm1 hypersurface} is now complete.\hfill $\Box$


\section{\textbf{The first eigenvalue of the focal submanifolds}}

At the beginning of this section, we should investigate the
multiplicity of the dimension $n-m_i$ as an eigenvalue of the focal
submanifold $M_i$ ($i=1,2$) of an isoparametric hypersurface with
$g$ distinct principal curvatures. For this purpose, we first
prepare the following lemma.

\begin{lem}\label{full immersion}
\emph{Both $M_1$ and $M_2$ are fully embedded in $S^{n+1}(1)$ if $g\geq 3$, namely, they cannot be embedded into
a hypersphere.
}
\end{lem}

\noindent
\textbf{\emph{Proof}}. We are mainly concerned with the proof for $M_1$; the other case is verbatim with obvious
changes on index ranges.

Suppose $M_1$ is not fully embedded in $S^{n+1}(1)$, then we can
find a point $q\in S^{n+1}(1)$ such that $\langle x, q \rangle=0$
for any $x\in M_1$. For any $p\in S^{n+1}(1)$, define the spherical
distance function $L_p: M_1\rightarrow \mathbb{R}$ by:
$$L_p(x)=\cos^{-1}\langle p,x\rangle.$$
Since $L_p$ is a Morse function on $M_1$ when $p\in S^{n+1}(1)-( M_1
\cup M_2)$ (\emph{cf.} \cite{CR}, p.285), we need only to deal with
the left two cases:

$(1)$ $p\in M_1$. Since the function $\langle x, p \rangle$ can
achieve $1$ at $x=p$, the point $q$ cannot lie in $M_1$.

$(2)$ $p\in M_2$. If $L_p$ is a constant, then from
each point $x\in M_1$, there exists one normal geodesic (normal to $M_1$ at $x$, normal to $M_2$ at $p$,
geodesic in $S^{n+1}(1)$), which connects $x$ and $p$.
Thus we can define a smooth map $f$ from the unit normal space of $M_2$ at $p$ to $M_1$ by:
\begin{eqnarray*}\label{Lp}
f: S(T^{\perp}_pM_2)\longrightarrow M_1\\
\xi \longmapsto x\quad
\end{eqnarray*}
where $x$ is the first intersection point of $M_1$ and the normal
geodesic starting from $p$ along the initial direction $\xi$ after
$\xi$ passes through the isoparametric hypersurface $M$. Under our
assumption, $f$ would be surjective. According to Sard Theorem, this
implies an inequality $m_2\geq \frac{g}{2}(m_1+m_2)-m_1$. Obviously,
this inequality holds true only when $g\leq 2$.

This completes the proof of Lemma \ref{full immersion}. \hfill $\Box$
\vspace{2mm}

\begin{rem}
The assumption $g\geq 3$ in Lemma \ref{full immersion} is essential.
For instance, for $g=2$, both the focal submanifolds of the
isoparametric hypersurface (generalized Clifford torus) are not
full, which are actually totally geodesic.
\end{rem}

As a direct result of Lemma \ref{full immersion}, the dimension $n-m_1$ (\emph{resp.} $n-m_2$) of $M_1$ is an eigenvalue
of $M_1$ (resp. $M_2$) with multiplicity at least $n+2$.
\vspace{4mm}

Now, we are ready to prove Theorem \ref{thm2 focal submanifold}.
\vspace{3mm}

\noindent \textbf{Theorem 1.3.}\, \emph{Let $M_1$ be the focal
submanifold of an isoparametric hypersurface with four distinct
principal curvatures in the unit sphere $S^{n+1}(1)$ with
codimension $m_1+1$. If $\dim M_1\geq \frac{2}{3}n+1$, then
$$\lambda_1(M_1)=\dim M_1$$ with multiplicity $n+2$.
A similar conclusion holds for $M_2$ under an analogous condition.}
\vspace{3mm}

\noindent
\textbf{\emph{Proof}.}
For sufficiently small $\varepsilon> 0$, set
$$M_1(\varepsilon):=S^{n+1}(1)- B_{\varepsilon} (M_2)=\bigcup_{\theta\in[0, \frac{\pi}{4}-\varepsilon]}M_{\theta}$$
where $B_{\varepsilon} (M_2)=\{x\in S^{n+1}(1)~|~dist (x, M_2)<\varepsilon\}$, $M_{\theta}$ is the isoparametric hypersurface with an oriented distance $\theta$ from $M_1$. Notice that the notation $M_{\theta}$ here is different from that we used before.

Given $\theta\in (0,\frac{\pi}{4}-\varepsilon]$, let
$\{e_{\alpha,i}~ \boldsymbol{|} ~ i=1,...,m_{\alpha},~
\alpha=1,..,4,~ e_{\alpha, i}\in E_{\alpha}\}$ be a local
orthonormal frame field on $M_{\theta}$ and $\xi$ be the unit normal
field of $M_{\theta}$ towards $M_1$. After a parallel translation
from any point $x\in M_{\theta}$ to a point $p=\phi_{\theta}(x)\in
M_1$, (where $\phi_{\theta}: M_{\theta}\rightarrow M_1$ is the focal
map, whose meaning is a little different from that in last section),
$\xi$ is still a unit normal vector at $p$, which we also denote by
$\xi$; $e_{1,i}$ $(i=1,...,m_1)$ turn to be normal vectors on $M_1$,
while the others are still tangent vectors on $M_1$, which we will
denote by $\{\widetilde{e}_{1,i}, \widetilde{e}_{2,i},
\widetilde{e}_{3,i}, \widetilde{e}_{4,i}\}$ determined by $x$.

For any $X\in T_xM_{\theta}$, we can decompose it as
$X=X_1+X_2+X_3+X_4\in E_1\oplus E_2\oplus E_3\oplus E_4$.
Identifying the principal distribution $E_{\alpha}(x)$
($\alpha=2,3,4$, $x\in M_{\theta}$) with its parallel translation at
$p=\phi_{\theta}(x)\in M_1$. The shape operator $A_{\xi}$ at $p$ is
given in terms of its eigenvectors $\widetilde{X}_{\alpha}$ (the
parallel translation of $X_{\alpha}, \alpha=2,3,4)$ by (\emph{cf.}
\cite{Mun})
\begin{eqnarray}\label{tilde X}
A_{\xi}\widetilde{X}_2&=&\cot(\theta_2-\theta_1)\widetilde{X}_2=\widetilde{X}_2,\nonumber\\
A_{\xi}\widetilde{X}_3&=&\cot(\theta_3-\theta_1)\widetilde{X}_3=0, \\
A_{\xi}\widetilde{X}_4&=&\cot(\theta_4-\theta_1)\widetilde{X}_4=-\widetilde{X}_4.\nonumber
\end{eqnarray}
Namely, $\widetilde{X}_2, \widetilde{X}_3, \widetilde{X}_4$ belong to the eigenspaces $E(1),E(0),E(-1)$ of $A_{\xi}$, respectively.

On the other hand, for a fixed $\theta$, define $\rho=\phi_{\theta}:
M_{\theta}\rightarrow M_1$. For any point $p\in M_1$, at a point
$x\in \rho^{-1}(p)$, we have a distribution $E_1\oplus E_2\oplus
E_3\oplus E_4$. Among them, the first one is projected to be $0$
under $\rho_{\ast}$; for the others, we have
\begin{eqnarray*}
\rho_{\ast}e_{\alpha, i}&=&\frac{\sin(\theta_{\alpha}-\theta)}{\sin\theta_{\alpha}}\widetilde{e}_{\alpha,i}
=\frac{\sin\frac{\alpha-1}{4}\pi}{\sin(\frac{\alpha-1}{4}\pi+\theta)}\widetilde{e}_{\alpha,i}\\
&:=&\widetilde{k}_{\alpha-1}\widetilde{e}_{\alpha,i},~~\qquad\qquad~i=1,...,m_{\alpha},~\alpha=2,3,4.
\end{eqnarray*}
Denote by $\{\theta_{\alpha,i}~|~\alpha=1,2,3,4,
i=1,...,m_{\alpha}\}$ the dual frame of $e_{\alpha,i}$.
We then conclude that (up to a sign)
\begin{equation}\label{volume Mtheta}
dM_{\theta}=\prod_{j=1}^{m_{\alpha}}\prod_{\alpha=2}^4\theta_{\alpha,j}\wedge
\prod_{i=1}^{m_1}\theta_{1,i}
=\frac{1}{\widetilde{k}_1^{m_2}\widetilde{k}_2^{m_1}\widetilde{k}_3^{m_2}}\rho^{\ast}(dM_1)\wedge
\prod_{i=1}^{m_1}\theta_{1,i}.
\end{equation}
Notice that here the submanifold $M_1$ may be non-orientable, but
the notation $dM_1$ still makes sense locally, up to a sign.

Let $h$ be the same function as in Section $2$. For sufficiently
small $\eta>0$, define $\widetilde{\psi}_{\eta}$ to be a nonnegative
smooth function on $[0, \frac{\pi}{2}-\eta]$ by
\begin{equation*}
\widetilde{\psi}_{\eta}(x):=\left\{ \begin{array}{ll}
1,\qquad\qquad x\in[0,\frac{\pi}{4}]\\
h(\frac{\frac{\pi}{2}-x}{\eta}),\quad x\in[\frac{\pi}{4}, \frac{\pi}{2}-\eta]
\end{array}\right.
\end{equation*}
Let $f_k$ $(k=0, 1, ...)$ be the $k$-th eigenfunctions on $M_1$ which are orthogonal to
each other with respect to the square integral inner product on $M_1$ and $L_{k+1} = Span\{f_0, f_1,..., f_k\}$.
Then any $\varphi\in L_{k+1}$ on $M_1$ can give rise to a function $\widetilde{\Phi}_{\varepsilon}: M_1(\varepsilon)\rightarrow \mathbb{R}$ by:
$$\widetilde{\Phi}_{\varepsilon}(x) = \widetilde{\psi}_{2\varepsilon}(2\theta)(\varphi\circ\rho)(x).$$
Evidently, similarly as last section,
$\widetilde{\Phi}_{\varepsilon}$ is a smooth function on
$M_1(\varepsilon)$ satisfying the Dirichlet boundary condition and
square integrable on $M_1(\varepsilon)$.

As in Section $2$, the calculation of $\|\nabla
\widetilde{\Phi}_{\varepsilon}\|^2_2$ is closely related to $|\nabla
\varphi(\rho)|^2$.
According to the decomposition (\ref{tilde X}), in the tangent space
of $M_1$ at $p$, we can decompose $\nabla \varphi$ as $\nabla
\varphi = Z_1+Z_2+Z_3\in E(1)\oplus E(0)\oplus E(-1)$. Thus we have
\begin{equation}\label{nabla varphi rho}
\left\{ \begin{aligned}
\quad|\nabla \varphi|_p^2~~~&=|Z_1|^2+|Z_2|^2+|Z_3|^2\\
|\nabla
\varphi(\rho)|_x^2&=\widetilde{k}^{2}_1|Z_1|^2+\widetilde{k}^{2}_2|Z_2|^2+\widetilde{k}^{2}_3|Z_3|^2
\end{aligned}\right.
\end{equation}

In the following, we will investigate the change of $|\nabla \varphi(\rho)|^2$ along with the point $x$ in the fiber sphere at $p$.
For this purpose, we recall

\noindent
\textbf{Lemma (see, for example, \cite{CCJ})}\,\,
{\itshape
Let $M^n$ be an isoparametric hypersurface in the unit sphere $S^{n+1}(1)$. Then the curvature
distributions are completely integrable. Their integral submanifolds corresponding to
$\cot\theta_j$ are totally geodesic in $M^n$ and have
constant sectional curvature $1 + \cot^2\theta_j$.
}

Denote by $S^{m_1}(\frac{1}{\sqrt{1+\cot^2\theta}}) \subset
M_{\theta}$ the fiber sphere at $p$. Clearly, for any pair of
antipodal points $x$, $x^{\prime}\in \rho^{-1}(p)=
S^{m_1}(\frac{1}{\sqrt{1+\cot^2\theta}})$, we have
$\xi(x^{\prime})=-\xi(x)$ by the parallel translations from $x$ and
$x^{\prime}$ to $p$, respectively. Denote by $E^{\prime}(1),
E^{\prime}(0), E^{\prime}(-1)$ the eigenspaces of
$A_{\xi(x^{\prime})}$ at $p$.
Then we can also decompose $\nabla \varphi$ as $\nabla \varphi = Z_3+Z_2+Z_1\in E^{\prime}(1)\oplus E^{\prime}(0)\oplus E^{\prime}(-1)$
with respect to $x^{\prime}$. In other words,
\begin{equation*}
\left\{\begin{array}{ll}
|\nabla \varphi(\rho)|^2_x = \widetilde{k}^{2}_1|Z_1|^2+\widetilde{k}^{2}_2|Z_2|^2+\widetilde{k}^{2}_3|Z_3|^2\\
|\nabla \varphi(\rho)|^2_{x^{\prime}} =
\widetilde{k}^{2}_3|Z_1|^2+\widetilde{k}^{2}_2|Z_2|^2+\widetilde{k}^{2}_1|Z_3|^2.
\end{array}\right.
\end{equation*}
Thus at the pair of two antipodal points $x$ and $x^{\prime}$, we
have
$$\frac{1}{2}\Big(|\nabla \varphi(\rho)|^2_x + |\nabla \varphi(\rho)|^2_{x^{\prime}}\Big)=\frac{\widetilde{k}^{2}_1+\widetilde{k}^{2}_3}{2}\Big(|Z_1|^2+|Z_3|^2\Big)+\widetilde{k}^{2}_2|Z_2|^2.$$

Set $\displaystyle \widetilde{K}:=\max\{\frac{\widetilde{k}^{2}_1+\widetilde{k}^{2}_3}{2}, \widetilde{k}^{2}_2\}$ for $\theta\in(0, \frac{\pi}{4}-\varepsilon].$ It is clear to see $\widetilde{K}=\frac{1}{\cos^22\theta}$ by the definition of $\widetilde{k}_{\alpha-1}$.
Since the assumption $3\dim M_1\geq 2n+3$ implies $m_2\geq 2$, which
guarantees that $\displaystyle \lim_{\varepsilon \rightarrow
0}\frac{1}{\varepsilon^2}\int_{\frac{\pi}{4}-2\varepsilon}^{\frac{\pi}{4}-\varepsilon}\cos^{m_2}2\theta
d\theta=0$. Then a similar discussion as in Section $2$ leads to
$$\lim_{\varepsilon \rightarrow
0}\int_{M_1(\varepsilon)}(\widetilde{\psi}^{\prime}_{2\varepsilon}(2\theta))^2
\varphi(\rho)^2 dM_1(\varepsilon)=0.$$
Hence
\begin{eqnarray}\label{nabla Phi prime}
\lim_{\varepsilon\rightarrow 0}\|\nabla
\widetilde{\Phi}_{\varepsilon}\|_2^2
&=&\lim_{\varepsilon\rightarrow
0}\int_{M_1(\varepsilon)}(\widetilde{\psi}_{2\varepsilon}(2\theta))^2
|\nabla (\varphi\circ \rho)|^2 dM_1(\varepsilon)\nonumber\\
&=&\int_{0}^{\frac{\pi}{4}}\Big(\int_{M_{\theta}}\frac{|\nabla
(\varphi\circ
\rho)|^2}{\widetilde{k}_1^{m_2}\widetilde{k}_2^{m_1}\widetilde{k}_3^{m_2}}\rho^{\ast}(dM_1)dS^{m_1}(\frac{1}{\sqrt{1+\cot^2\theta}})\Big)d\theta\\
&\leq&\int_{0}^{\frac{\pi}{4}}\Big(\int_{M_{\theta}}|\nabla
\varphi|^2\rho^{\ast}(dM_1)dS^{m_1}(\frac{1}{\sqrt{1+\cot^2\theta}})\Big)\cdot \frac{\widetilde{K}}{\widetilde{k}_1^{m_2}\widetilde{k}_2^{m_1}\widetilde{k}_3^{m_2}}d\theta\nonumber\\
&=&\int_{0}^{\frac{\pi}{4}}\Big(\int_{M_1}|\nabla \varphi|^2 dM_1
\Big)\cdot Vol(S^{m_1}(\frac{1}{\sqrt{1+\cot^2\theta}}))\cdot \frac{\widetilde{K}}{\widetilde{k}_1^{m_2}\widetilde{k}_2^{m_1}\widetilde{k}_3^{m_2}}d\theta\nonumber\\
&=&\|\nabla \varphi\|_2^2\cdot\frac{C_{m_1}}{2^{m_1+1}}\int_0^{\frac{\pi}{2}}\sin^{m_1}\theta\cos^{m_2-2}\theta~d\theta\nonumber\\
&=&\|\nabla \varphi\|_2^2\cdot\frac{C_{m_1}}{2^{m_1+2}}\cdot B(\frac{m_1+1}{2}, \frac{m_2-1}{2}).\nonumber
\end{eqnarray}
where $Vol(S^{m_1}(\frac{1}{\sqrt{1+\cot^2\theta}}))=C_{m_1}\cdot\sin^{m_1}\theta$, $C_{m_1}$ is the volume of $S^{m_1}(1)$.
Besides, with a simple calculation, we get
\begin{eqnarray}\label{Phi prime}
\lim_{\varepsilon\rightarrow 0}\|\widetilde{\Phi}_{\varepsilon}\|_2^2
&=&\int_0^{\frac{\pi}{4}}\frac{1}{\widetilde{k}_1^{m_2}\widetilde{k}_2^{m_1}\widetilde{k}_3^{m_2}} \int_{M_1}\int_{S^{m_1}(\frac{1}{\sqrt{1+\cot^2\theta}})} \varphi(\rho)^2 dS^{m_1}dM_1d\theta\nonumber\\
&=& \|\varphi\|_2^2\cdot\int_0^{\frac{\pi}{4}}\frac{1}{\widetilde{k}_1^{m_2}\widetilde{k}_2^{m_1}\widetilde{k}_3^{m_2}} Vol(S^{m_1})~d\theta \\
&=&\|\varphi\|_2^2\cdot \frac{C_{m_1}}{2^{m_1+2}}\cdot B(\frac{m_1+1}{2}, \frac{m_2+1}{2}).\nonumber
\end{eqnarray}
Consequently, combing with (\ref{nabla Phi prime}) and (\ref{Phi prime}), we arrive at
\begin{equation*}
\lim_{\varepsilon\rightarrow 0}\frac{\|\nabla \widetilde{\Phi}_{\varepsilon}\|_2^2}{\|\widetilde{\Phi}_{\varepsilon}\|_2^2}
\leq \frac{\|\nabla \varphi\|_2^2}{\|\varphi\|_2^2}\cdot \frac{B(\frac{m_1+1}{2}, \frac{m_2-1}{2})}{B(\frac{m_1+1}{2}, \frac{m_2+1}{2})}
=\frac{\|\nabla \varphi\|_2^2}{\|\varphi\|_2^2}\cdot \frac{m_1+m_2}{m_2-1}.
\end{equation*}
A similar argument as in Section $2$ leads us to
\begin{equation}\label{frac Phi prime}
\lambda_k(S^{n+1}(1))\leq \lambda_k(M_1) \frac{m_1+m_2}{m_2-1}.
\end{equation}
This inequality connects the eigenvalues of $S^{n+1}(1)$ and that of the focal submanifold $M_1$ in a concise manner.
It contains rich information. Now we take $k=n+3$. The inequality (\ref{frac Phi prime}) turns to
$$\frac{2(n+2)(m_2-1)}{m_1+m_2}\leq \lambda_{n+3}(M_1).$$
Based on this inequality, in order to complete the proof of Theorem \ref{thm2 focal submanifold}, we just need
to establish the following inequality
\begin{equation}\label{lambda1 of M1}
\dim M_1=m_1+2m_2< \frac{2(n+2)(m_2-1)}{m_1+m_2}.
\end{equation}
Due to the relation $n=2(m_1+m_2)$, we get a sufficient
condition on the positive integers $m_1$, $m_2$ which is almost optimal for the inequality (\ref{lambda1 of M1}) to hold:
$$m_2\geq\frac{1}{2}(m_1+3).$$
At last, combing with Lemma \ref{full immersion}, we can conclude that
$$\lambda_1(M_1)=\dim M_1=m_1+2m_2,~ with~ multiplicity~n+2, \quad provided~ m_2\geq\frac{1}{2}(m_1+3)$$
as we required.
\hfill $\Box$

\vspace{4mm}

\begin{rem}
When $g=1$, the focal submanifolds are just two points. When $g=2$,
as is well known, the isoparametric hypersurface in $S^{n+1}(1)$ is
isometric to the generalized Clifford torus
$S^p(\sqrt{\frac{p}{n}}~)\times S^{q}(\sqrt{\frac{q}{n}}~)$
$(p+q=n)$. The focal submanifolds are isometric to $S^p(1)$ and
$S^q(1)$. Clearly, their first eigenvalues are their dimensions.
When $g=3$, E. Cartan asserted that $m_1=m_2=1, 2, 4~or~ 8$. The
focal submanifolds in the unit sphere $S^4(1)$, $S^7(1)$,
$S^{13}(1)$ and $S^{25}(1)$ are the Veronese embedding of
$\mathbb{R}P^2$, $\mathbb{C}P^2$, $\mathbb{H}P^2$ and
$\mathbb{O}P^2$, respectively. For this $\mathbb{R}P^2$ minimally
embedded in $S^4(1)$, its induced metric differs the standard metric
of constant Gaussian curvature $K=1$ by a constant factor such that
$K=\frac{1}{3}$, thus $\lambda_1(\mathbb{R}P^2)=2.$ As for these
$\mathbb{C}P^2$, $\mathbb{H}P^2$ and $\mathbb{O}P^2$, they are
minimally embedded in the unit spheres $S^7(1)$, $S^{13}(1)$ and
$S^{25}(1)$, respectively, while the induced metric differs the
symmetric space metric by a constant factor such that
$\frac{1}{3}\leq Sec \leq \frac{4}{3}$. By \cite{Str} and
\cite{Mas}, the first eigenvalues of the focal submanifolds
$\mathbb{C}P^2$, $\mathbb{H}P^2$ and $\mathbb{O}P^2$ are equal to
their dimensions, respectively.

Therefore, for $g=2, 3$,
$$\lambda_1(M_i)=\dim M_i,\quad i=1,2.$$\hfill $\Box$
\end{rem}

We conclude this paper with a proof of Proposition \ref{M1 for other}.
\vspace{3mm}

\noindent
\textbf{Proposition 1.1.}\emph{Let $M_2$ be the focal submanifold of OT-FKM-type defined before with $(m_1, m_2)=(1, k)$. The following equality is valid
\begin{equation*}
\lambda_1(M_2)=\min\{4, 2+k\}.
\end{equation*}}
\vspace{3mm}

\noindent \textbf{\emph{Proof.}} When $m_1=1$, $m_2=k$, the
OT-FKM-type polynomial can be written as
\begin{eqnarray*}
&&\qquad \quad F: \mathbb{R}^{2k+4}\longrightarrow \mathbb{R}\\
&&F(x)=|x|^4-2(\langle P_0x,x \rangle^2+\langle P_1x,x \rangle^2).
\end{eqnarray*}
By orthogonal transformations, we can always choose $P_0$ and $P_1$
to be
\begin{eqnarray}\label{Clifford A_j}
&&P_0= \left(\begin{array}{c|c}
I & 0 \\
\hline 0 & -I
\end{array}\right),\quad
P_1= \left(\begin{array}{c|c}
0 & I \\
\hline I & 0
\end{array}\right).
\end{eqnarray}
Writing any point $x\in S^{2k+3}(1)$ as $x=(z,w)\in
\mathbb{R}^{k+2}\times\mathbb{R}^{k+2}$, the focal submanifold
$M_2=f^{-1}(-1)$ ($f=F|_{S^{2k+3}(1)}$) can be characterized as
\begin{equation*}
M_2^{k+2}=\{(z,w)\in S^{2k+3}(1)~|~z \sslash w\}.
\end{equation*}

Define a map
\begin{eqnarray*}
\Psi: S^1(1)\times S^{k+1}(1)&\longrightarrow& M_2^{k+2}\subset \mathbb{R}^{2k+4}\\
e^{i\theta}, x=(x_1,...,x_{k+2})&\mapsto& (e^{i\theta}x_1,...,e^{i\theta}x_{k+2}).
\end{eqnarray*}
It satisfies $\Psi(\theta+\pi, -x)=\Psi(\theta, x)$. In this way, we
can identify isometrically $M_2$ with the metric induced from
$S^{2k+3}(1)$ as
$$M_2^{k+2}\cong S^1(1)\times S^{k+1}(1)\Big/(\theta,
x)\sim(\theta+\pi,-x).$$

The eigenfunctions of $M_2$ are those products of eigenfunctions
from $S^1(1)$ and $S^{k+1}(1)$ which take the same values at
$(\theta, x)$ and $(\theta+\pi,-x)$.
Hence $\lambda_1(M_2^{k+2})=\min\{4, k+2\}$, as we claimed.



\vspace{8mm}


\begin{ack}
The authors would like to thank Professors T. Cecil and C. K. Peng for their helpful comments on OT-FKM isoparametric foliation
and spectrum of Laplacian, respectively. We also express our gratitude to Professors Y. B. Shen and W. P. Zhang for valuable discussions and
Professors Q. M. Cheng and Y. Ohnita for their interests.
\end{ack}

\end{document}